\documentclass[numreferences]{kluwer}
\usepackage{amssymb}
\usepackage{latexsym}

\newtheorem{thm}{Theorem}[section] 
\newtheorem{dfn}[thm]{Definition}
\newtheorem{rmk}[thm]{Remark}

\newtheorem{prop}[thm]{Proposition}
\newtheorem{lem}[thm]{Lemma}

\newcommand{\Pf}{{\em Proof}. }
\newcommand{\EPf}
{%
\mbox{}%
\nolinebreak%
\hfill%
\rule{2mm}{2mm}%
\medbreak%
\par%
}

\newcommand{\C}{\mathbb C} 
\newcommand{\R}{\mathbb R}

\newcommand{\g}{{\cal G}{}}

\newcommand{\s}{{\cal S}{}} 

\renewcommand{\d}{{\cal D}{}}

\renewcommand{\o}{{\cal O}{}}

\renewcommand{\Re}{\mbox{Re}}
\newcommand{\sm}{\star_{\nu}^{M}}
\newcommand{\sw}{\star_{q}^{W}}
\newcommand{\sn}{\star_{\nu}}
\newcommand{\skq}{\star_{q}^{(k)}}

\newcommand{\ronu}{\rho_{\nu}}
\newcommand{\ronuhat}{{\hat{\rho}}_{\nu}}
\newcommand{\LT}{{\cal L}{}}
\newcommand{\FT}{{\cal F}{}}
\newcommand{\ZT}{{\cal Z}_\nu{}}
\newcommand{\taukq}{\tau^{(k)}_q{}}
\newcommand{\Tkq}{T^{(k)}_q{}}

\def\Der{{\cal	D}er}

\newcommand{\sab}{{\cal S}_\alpha^\beta{}}
\newcommand{\sba}{{\cal S}_\beta^\alpha{}}
\newcommand{\sA}{{\cal S}^{\sigma(\alpha_{1},\alpha_{2})}_{(\alpha_{1},\alpha_{2})}(2){}}
\newcommand{\Ekab}{E^{(k)}_{(\alpha_{1},\alpha_{2})}{}}

\title{Convergent star product algebras on ``$ax+b$''}
\author{{\bf Pierre Bieliavsky}\\
Universit\'e Libre de Bruxelles, Belgium\\
e-mail: pbiel@ulb.ac.be\\
{\bf Yoshiaki Maeda}\\
Keio University, Japan\\
e-mail: maeda@math.keio.ac.jp}
\begin{document}
\maketitle
\section{Introduction}
The notion of convergent star product is generally understood as the 
data of a one parameter family $\{E_t\}_{t\in I}\subset C^\infty(M)$ of 
function algebras on a Poisson manifold $(M,\{\, , \,\})$.  On each of 
them
one is given an associative algebra 
structure $\star_t$ which respect to which the function space $E_t$ is 
closed. The family of products $\{\star_t\}$ should moreover define in some 
sense a deformation of the commutative pointwise product of functions in 
the direction of the Poisson structure $\{\, , \,\}$.

Stable function algebras for the Weyl-Moyal product have been studied in 
various contexts. For instance, see:
\begin{enumerate}
\item[-] \cite{M} for such a study in the framework of tempered 
distributions on a symplectic vector space;
\item[-] \cite{R} for a $C^\star$-algebraic study on 
$\R^d$-manifolds;
\item[-] \cite{OMMY} for non-tempered stable function spaces on $\C^n$.
\end{enumerate}
A special feature of the Weyl-Moyal star product---independently of 
the functional framework--- is its maximal invariance under the group 
of affine transformations with respect to a {\sl flat} affine 
connection. This can be rephrased by saying that the Weyl-Moyal 
quantization is universal with respect to actions of $\R^d$  
\cite{R}. A natural question is then the one of defining universal 
(convergent) deformations for non-Abelian Lie group actions. 
In the formal framework, this has been investigated in \cite{GZ}.
In, \cite{BB,H}, such formulae in the case of $ax+b$ have been studied within the context of
Wigner formalism and signal analysis. However, the question of defining
an adapted functional framework has not been investigated. In 
\cite{B,BM}, one finds a functional analytic study for solvable Lie 
group actions and symmetric spaces in the tempered distributions and/or $C^\star$-context.

It therefore appears quite challenging to investigate the problem of 
defining non-tempered (e.g. exponential growth) function spaces on 
such a non-Abelian Lie group which are stable under some left-invariant 
(convergent) star product. In other words, studying a situation where 
non-temperedness and non-linear invariance mix. This is what is done 
in this paper for the particular case of the group $ax+b$. More 
precisely, we first start by giving a construction of a left-invariant 
star product on the (symplectic) group manifold underlying the Lie 
group $ax+b$.  This star product is obtained via an equivalence 
transformation $T$ performed on Moyal's product (which is not 
left-invariant). The equivalence $T$ involves two ingredients: a 
partial Laplace transform and a family 
$\{\phi_{\nu,\gamma}:\C\to\C\}_{\nu,\gamma\in\C}$ of holomorphic maps.
For special values of of the parameters $\nu$ and $\gamma$ one refinds 
the functional calculus studied in \cite{B,BM}. But for other values, 
one can define stable function algebras constituted by type-S functions.
Using a holomorphic presentation of 
these spaces, one gets non-commutative algebra structures on spaces of
entire functions on $\C^2$. Such a space contains functions of 
exponential growth, one therefore has a non-tempered invariant 
calculus.

\noindent{\bf Acknowledgment} We thank Daniel Sternheimer for having 
communicated important references.
\section{Formal star products on ``$ax+b$''}\label{FORMAL}
In this section, we briefly recall results appearing in \cite{B,BM}.
Let $\g$ denote the Lie algebra of the group of affine 
transformations of the real line. Formally, one has 
$\g=\mbox{span}_{\R}\{A,E\}$ with table $[A,E]=2E$. Consider 
the linear map $\lambda:\g\to C^\infty(\R^2):X\mapsto\lambda_X$
defined by $\lambda_A(a,l)=2l;\quad\lambda_E(a,l)=e^{-2a}$, where 
$\R^2=\{(a,l)\}$. One then checks that the map $\lambda$ is a homomorphism 
of Lie algebras when $C^\infty(\R^2)$ is endowed with the symplectic Poisson 
bracket $\{\, , \,\}:=\partial_a\wedge\partial_l$. Moreover, if $\sm$ 
denotes the formal Moyal star product on $C^\infty(\R^2)[[\nu]]$ 
(i.e 
$u\sm v=u\,\exp{(\nu\stackrel{\leftarrow}{\partial_a}\wedge
\stackrel{\rightarrow}{\partial_l})}v\qquad u,v\in C^\infty(\R^2)[[\nu]]$),
one has $\left[\lambda_A,\lambda_E\right]_\nu=2\nu\{\lambda_A,\lambda_E\}$
(where $\left[u,v\right]_\nu:=u\sm v-v\sm u$). In particular, the formula
$$
\ronu(X)u:=\frac{1}{2\nu}\left[\lambda_X,u\right]_\nu\qquad X\in\g,u\in C^\infty(\R^2)[[\nu]]
$$
defines a homomorphism of Lie algebras 
$$
\ronu:\g\to\Der(C^\infty(\R^2)[[\nu]],\sm).
$$
Explicitly, one has 
$\ronu(A)u=-\partial_au;\quad\ronu(E)u=-\frac{e^{-2a}}{\nu}\sinh(\nu\partial_l)u$.
Intertwining the representation $\ronu$ by a transformation of the type
$$
\LT(u)(a,z):=\int_\R e^{-zl}u(a,l)\,dl,
$$
one gets
$$
\begin{array}{ccc}
\ronuhat(A)\LT(u):=\LT(\ronu(A)u) & = & -\partial_a\LT(u);\\
\ronuhat(E)\LT(u):=\LT(\ronu(E)u) & = & -\frac{e^{-2a}}{\nu}\sinh(\nu 
z)\LT(u),
\end{array}
$$
where we assumed $u(a,\pm\infty)=0$. Now, set formally
$$
\ZT(u)(a,z):=\int_\R \exp{\left(\gamma\frac{1}{\nu}\sinh(\nu z)l\right)}\,u(a,l)\,dl,
$$
and 
$$
f\bullet_\nu g:=\ZT(\ZT^{-1}f.\ZT^{-1}g)\qquad(\gamma\in\C_{0}).
$$
\begin{prop}
For all $X\in\g$, $\ronuhat(X)$ is a derivation of the commutative 
product $\bullet_\nu$.
\end{prop}
In other words, the associative formal product $u\sn v:=T^{-1}(Tu\sm Tv)$ 
where $T=\LT^{-1}\ZT$ is invariant under the infinitesimal action 
$\g\to\Gamma T(\R^2):X\mapsto{}^\sharp(d\lambda_X)$ where ${}^\sharp(d\lambda_X)$ 
denotes the Hamiltonian vector field associated (via the symplectic 
structure) to the function $\lambda_X$. This action of $\g$ turns out to 
exponentiate as a global simply transitive symplectic action of the group $G=``ax+b"$ on 
$\R^2$, providing an identification of the group manifold 
underlying $G$ with $\R^{2}$.

The integral form of the transformation $\ZT$ allows to define specific 
functions algebras on $\R^2$ (as opposed to power series  algebras) stable 
under the product $\sn$ where the formal Moyal product $\sm$ is replaced 
by its ``convergent" version: the Weyl product. The case where $\nu\in 
i\R,\quad\gamma=1,\quad z\in i\R$ has been studied in \cite{BM,B}. In 
what follows, we are concerned with the general case where $\nu\in 
i\R,\quad\gamma\in U(1),\quad z\in \C$. We end up this section by observing 
that the intertwiner $T=\LT^{-1}\ZT$ can be expressed as 
$$
T=\LT^{-1}\circ(\phi_{\nu,\gamma})^*\circ\LT
$$
where we set $\phi_{\nu,\gamma}(z):=\frac{\gamma}{\nu}\sinh(\nu z)$. The 
map $\phi_{\nu,\gamma}$ will be referred in the sequel as the ``twisting 
map" (see Section \ref{TM}).
\begin{rmk}
{\rm
An alternative simple way for obtaining an explicit formula of an 
invariant star product on $ax+b$ is based on the following observation
\cite{BB,H}. The symplectic group manifold underlying $ax+b$ can be 
seen as an open coadjoint orbit $\cal O$ in $\g^{\star}$. Quantizing 
the 
Poisson manifold $\g^{\star}$ via the universal enveloping algebra 
product and  then restricting to $\cal O$ yields an invariant star 
product on $\cal O$ hence a left-invariant one on $ax+b$. It is classical 
that an oscillatory integral formula for this product can be written 
down in terms of the Campbell-Baker-Haussdorff function (see 
\cite{BB} for explicit computation: Formulae 5.8 and 5.12). Our product
$\star_{\nu}$ described above is different from the universal 
enveloping algebra product. Indeed, their invariance diffeomorphism 
groups do not coincide \cite{B}.
}
\end{rmk}
\section{Fundamental spaces of type ${\cal S}$}\label{TYPE-S}
In this section, we follow Chapter IV of I.M.Guelfand's book \cite{G}.
We will denote by $\o(\C^m)$ the space of holomorphic (entire) functions 
on $\C^m$. 
\begin{dfn} Let $\alpha,\beta\in\R^m$.
The fundamental space $\sab(m)$ is defined as the space of holomorphic 
functions $\varphi\in\o(\C^m)$ such that there exists $a,b\in(\R^+)^m$ 
and $C>0$ with
$$
|\varphi(x+iy)|\leq 
C\,\exp\left(-a|x|^{\frac{1}{\alpha}}+b|y|^{\frac{1}{1-\beta}}\right),
$$
where we adopt the usual notations~:
$a|x|^e=\sum_ja_j|x_j|^{e_j}\quad(a,x,e\in\R^m); 
\,\frac{1}{\alpha}=(\frac{1}{\alpha_1},...,\frac{1}{\alpha_m});\, 
1-\beta=(1-\beta_1,...,1-\beta_m)$.
\end{dfn}
Every element $\varphi\in \sab(m)$ is entirely determined by its restriction 
the ``real axis" $\varphi(x)\quad x\in\R^m$. We will often identify 
the space $\sab(m)$ with the subspace $\left(\sab(m)\right)_x$ of 
$C^\infty(\R^m)$ constituted by the restrictions. In order to consider 
only non-trivial spaces, we will assume 
$\alpha+\beta\geq1;\quad\alpha>0;\quad\beta>0$. We will denote by 
$\FT(u)(\xi)$ the Fourier transform of the function $u\in L^1(\R^m)$~:
$$
\FT(u)(\xi):=\int_{\R^m}e^{i\xi.x}u(x)\, dx ,
$$
where $\xi.x$ denotes the canonical dot product on $\R^m$. For even $m=2n$, 
we will denote by $J$ the endomorphism of $\R^{2n}$ defined by the matrix
$$
[J]:=\left(
\begin{array}{cc}
0 & I_n\\
-I_n & 0
\end{array}\right)
$$
where $I_n$ is the $n\times n$ identity matrix. We denote by $\omega$ the 
bilinear symplectic structure on $\R^{2n}$ defined by $\omega(x,y):=x.Jy$.
\begin{dfn}
We define the symplectic Fourier transform of the function $u\in 
L^1(\R^{2n})$  as 
$$
S\FT(u)(y):=\int_{\R^{2n}}e^{i\omega(x,y)} u(x)\, dx \qquad (y\in\R^{2n}).
$$
\end{dfn}
Equivalently, one has $S\FT=J^*\circ\FT$, which yields (see \cite{G})
\begin{lem}
One has $$S\FT(\sab(2n))=\s^{\sigma(\alpha)}_{\sigma(\beta)}(2n),$$ where for 
$\alpha=(\alpha_1,\alpha_{2})\in\R^{{2n}}\,(\alpha_{i}\in\R^n)$, we set
$\sigma(\alpha):=(\alpha_2,\alpha_1)$.
\end{lem}
\begin{dfn}
For $u,v\in L^1(\R^{2n})$, one defines their twisted convolution by 
$$
u\times_qv(x):=\int_{\R^{2n}}e^{iq\omega(x,y)}u(y)\, 
v(x-y)\,dy\qquad(q\in\R_0).
$$
\end{dfn}
\begin{lem}\label{MULT}
One has $$\sab(m).\s^{\beta'}_{\alpha'}(m)\subset
\s^{\mbox{\rm max}(\beta,\beta')}_{\mbox{\rm min}(\alpha,\alpha')}(m);$$ where we set
$\mbox{\rm max}(\beta,\beta'):=(\mbox{\rm 
max}(\beta_1,\beta'_1),...,\mbox{\rm max}(\beta_m,\beta'_m))$.
\end{lem}
\Pf
For simplicity, we assume $m=1$. Let $\varphi\in\sab$ and $\varphi'\in\s^{\beta'}_{\alpha'}$.
Then, for $z=x+iy$, one has
{\small{
\begin{eqnarray*}
|\varphi\varphi'(z)|\leq CC'\exp\left(
-a|x|^{\frac{1}{\alpha}}+b|y|^{\frac{1}{1-\beta}}-a'|x|^{\frac{1}{\alpha'}}+
b'|y|^{\frac{1}{1-\beta'}}\right)
\end{eqnarray*}
}}
which is lower than 
{\small{
$C"\exp\left(
-a"|x|^{\mbox{\rm max}(\frac{1}{\alpha},\frac{1}{\alpha'})}
+b"|y|^{\mbox{\rm 
max}(\frac{1}{1-\beta},\frac{1}{1-\beta'})}\right)$}} for 
some $C",a",b"$. But, 
$\mbox{max}(\frac{1}{\alpha},\frac{1}{\alpha'})=\frac{1}{\mbox{\rm min}({\alpha},{\alpha'})}$ and 
$\mbox{max}(\frac{1}{1-\beta},\frac{1}{1-\beta'})=\frac{1}{1-\mbox{max}(\beta,\beta')}$.
\EPf
\begin{lem}
Let $u,v\in\sab(2n)$. Then
$u\times_qv\in\s^{\sigma{\alpha}}_{\sigma{\beta}}(2n)$.
\end{lem}
\Pf
Changing the variables following $y\mapsto-y$, one gets 
$u\times_qv=\left[d_q^* S\FT(\tilde{u}\alpha_xv)\right]^{\tilde{}}$ where 
$\tilde{u}(x):=u(-x)$, $(\alpha_xv)(y):=v(y+x)$ and where $d_q$ denotes the 
dilation in $\R^m$~: $x\mapsto qx \quad(q\in\R)$. Hence $u\times_qv\in 
S\FT\left(\s^{\mbox{\rm max}(\beta,\beta)}_{\mbox{\rm 
min}({\alpha},{\alpha})}(2n)\right)=\s^{\sigma(\alpha)}_{\sigma(\beta)}(2n)$.
\EPf
\begin{dfn}(see e.g. \cite{M})
The Weyl product between $u$ and $v$ in $L^{1}(\R^{2n})$ is defined by
$$
u\sw v:=S\FT\left[S\FT(u)\times_qS\FT(v)\right].
$$
\end{dfn}
\begin{prop}\label{WEYL}\cite{OMMY}
Let $u,v\in\sab(2n)$, Then $u\sw v\in\s_{\sigma(\beta)}^{\sigma(\alpha)}(2n)$. 
In particular, the space $\s_\alpha^{\sigma(\alpha)}(2n)$ is stable under 
the Weyl product.
\end{prop}
\begin{rmk}
{\rm
The space $\s_\alpha^{\sigma(\alpha)}(2n)$ is stable under the pointwise 
multiplication as well.
}
\end{rmk}
\section{Laplace Transformation}
In this section we follow L. Schwartz' book \cite{S}. We adopt the 
following notations. We denote by $\d$ the space of compactly supported 
smooth functions on $\R$ endowed with the topology of test 
functions. We denote by $\d'$ the space of distributions on $\R$.
Also, if $\Omega$ is an open domain in $\C$, we set $\o(\Omega)$ for the 
space of holomorphic functions on $\Omega$.
\begin{dfn}Let $\Gamma$ be an open interval in $\R$.
The fundamental space $\s'_\Gamma$ is defined as the space of 
distributions $T\in\d'$ such that for all $\xi\in\Gamma$, the distribution 
$\exp(-\xi.x)T_x$ is tempered. We denote by $\o_\Gamma$ the space of 
holomorphic functions $F\in\o(\Gamma+i\R)$ such that for all compact set 
$K\subset\Gamma$, the restriction $F|_{K+i\R}(\xi+i\eta)$ is bounded by a 
polynomial in $\eta$.
\end{dfn}
\begin{prop}
One defines the Laplace transform of an element $T\in\s'_\Gamma$ as the Fourier 
transform of $e^{-\xi.x}T_x$
$$
\LT(T)(\xi+i\eta):=\left(\FT_x\left[\exp(-\xi.x)T_x\right]\right)(\eta).
$$
Then, setting $z=\xi+i\eta\in\Gamma+i\R$, one has a linear isomorphism
$$
\LT:\s'_\Gamma\to\o_\Gamma.
$$
\end{prop}
\begin{rmk}
{\rm
Provided the following integrals make sense, one has
$$
\LT(T)(z)=\int_\R e^{-zx}T_x\,dx\mbox{ and }
\LT^{-1}F(x)=\int_{c+i\R}e^{zx}F(z)dz,
$$
where $c$ is any element of $\Gamma$. Indeed, if $F=\LT T$ one has for all 
$c\in\Gamma$~:
{\small{
\begin{eqnarray*}
\int_{c+i\R}e^{zx}F(z)dz=e^{xc}\int_\R F(c+it) e^{ixt}dt=\\
=
e^{xc}\int_\R\int_\R e^{it(x-y)}\left(T_ye^{-yc}\right)dy\,dt
=
e^{xc}\delta_x\left(T_ye^{-yc}\right)=T_x.
\end{eqnarray*}
}}
}
\end{rmk}
Now we assume $m=1$ and set $\sab(1)=:\sab$.
\begin{lem}
For all $\xi\in\R$ the function $e^{-\xi x}$ is a multiplier in  
$\sab$ ($\alpha<1$).
\end{lem}
\Pf
Let $\varphi\in\sab$. Then the function $f(z):=e^{-\xi z}\varphi(z)$ is 
entire and one has $$|f(z)|\leq C\exp\left(-a|x|^{\frac{1}{\alpha}}
+b|y|^{\frac{1}{1-\beta}}\right)e^{|\xi||x|}$$ which is lower than 
$C'\exp\left(-a'|x|^{\frac{1}{\alpha}}
+b|y|^{\frac{1}{1-\beta}}\right)$ as soon as $\alpha<1$.
\EPf
In particular, one has $\sab\subset\s'_\R=\cap_\Gamma\s'_\Gamma$ as soon as 
$\alpha<1$.
\begin{prop}
The Laplace transformation yields a linear isomorphism
$$
\LT:\left(\sab\right)_x\to J^*\sba
$$
where $J:\C\to\C$ denotes the multiplication by $i=\sqrt{-1}$.
\end{prop}
\Pf
Let $\varphi\in\sab$. Then 
$(J^*\LT(\varphi))(x+iy)=\int_\R e^{ixt-yt}\varphi(t)\,dt$. Hence 
$(J^*\LT(\varphi))|_\R(x)=(\FT(\varphi))(x)\in\sba$ (\cite{G}). Thus 
$\LT(\varphi)\in J^*\sba$. Now let $f\in J^*\sba$. One has 
$|f(x+iy)|\leq C\exp\left(-a|y|^{\frac{1}{\alpha}}
+b|x|^{\frac{1}{1-\beta}}\right)$ which guarantees that 
$f\in\o_\R$. Therefore
$\LT^{-1}f(x)=e^{xc}\int_\R e^{ixt}f(c+it)\,dt\quad(c\in\R)$. Choosing 
$c=0$, one gets $\LT^{-1}f=\FT J^*f\in\sab$. One therefore has an 
isomorphism $\LT^{-1}:J^*\sba\to\sab$.
\EPf
\section{Twisting maps}\label{TM}
\begin{dfn} Let $q\in\R$ and $\theta\in[0,2\pi[$. We define the twisting 
map $\phi_{q,\theta}:\C\to\C$ by 
$$
\left\{
\begin{array}{ccc}
\phi_{q,\theta}(z) & = & \frac{e^{i\theta}}{q}\sinh(iqz) \mbox{ if } q\neq 
0\\
\phi_{0,\theta}(z) & = & z.
\end{array}
\right.
$$
\end{dfn}
\begin{lem}\label{TWIST}
The twisting map $\phi_{q,0}\quad(q\neq 0)$ establishes a biholomorphic 
diffeomorphism 
$$
\phi_{q,0}:S_q:=\{z\in\C|\,|\Re{(z)}|<\frac{\pi}{2q}\}\to
\C-\{\pm i[\frac{1}{q},\infty[\}.
$$
\end{lem}
\Pf
In coordinates $z=x+iy$, the twisting map is
$$\phi_{q,0}(x,y)=\left(\frac{1}{q}\sinh(qy)\cos(qx),\frac{1}{q}\cosh(qy)\sin(qx)\right).$$
In particular, the imaginary ax\-is $x=0$ is sent onto the real ax\-is $y=0$, 
while the image of the vertical line $x=\pm\frac{\pi}{2q}$ under $\phi_{q,0}$ 
is the half imaginary line $\pm i[\frac{1}{q},\infty[$. At last, for 
$c\in]0,\frac{\pi}{2}[$ the image of the vertical line $x=\pm\frac{c}{q}$ 
is the branch of hyperbola $\{-\left(\frac{qx}{\cos 
c}\right)^2+\left(\frac{qy}{\sin c}\right)^2=1\}\cap\{\pm y>0\}$. For 
$\epsilon\in]0,\frac{\pi}{2}[$, one therefore gets a holomorphic 
diffeomorphism between the strip 
$]-\frac{\epsilon}{q},\frac{\epsilon}{q}[+i\R$ and the region 
$\{-\left(\frac{qx}{\cos 
\epsilon}\right)^2+\left(\frac{qy}{\sin \epsilon}\right)^2<1\}$.
\EPf
As explained in Section~\ref{FORMAL}, we are interested in considering 
transformations of the type $\LT^{-1}\circ\phi_{q,\theta}^*\circ\LT$ or 
$\LT^{-1}\circ(\phi_{q,\theta}^{-1})^*\circ\LT$. Let $F$ be a function 
defined on some domain $\Omega$ of $\C$. In order to define 
$\LT^{-1}(\phi_{q,\theta}^{-1})^*F$, we want 
$(\phi_{q,\theta}^{-1})^*F\in\o_\Gamma$. In particular we want $(\phi_{q,\theta}^{-1})^*F$
to be defined on vertical lines. This imposes to $\theta$ to be an integral 
multiple of $\frac{\pi}{2}$.
\begin{prop}\label{ISOTWIST}
Let $\alpha,\beta\in]0,1[$ be such that $\alpha+\beta\geq1$.
\begin{enumerate}
\item[(i)] For all open interval $I$ of positive numbers $(I\subset\R_0^+)$, one has 
the injection $(\phi_{q,0}^{-1})^*:J^*\sba\to\o_I\cap\o_{-I}$.
\item[(ii)] For all open interval 
$I\subset]-\frac{\pi}{2q},\frac{\pi}{2q}[$, one has the injection
$(\phi_{q,\frac{\pi}{2}}^{-1})^*:J^*\sba\to\o_I$.
\end{enumerate}
\end{prop}
\Pf
Let $F\in J^*\sba$ and set for simplicity $\phi=\sinh(iz)$. Let 
$I\subset]0,\infty[$ and $K$ be a compact set in $I$. Then $\phi^{-1}$ is 
well defined on the strip $K+i\R$ and $\phi^{-1}(K+i\R)\subset]-\frac{\pi}{2},
\frac{\pi}{2}[+i\R$ (cf. Lemma~\ref{TWIST}). Now, consider the integral
$I_K:=\int_{K+i\R}|(\phi^{-1})^*F(z)|\,|dz\wedge\overline{dz}|$. If 
$I_K<\infty$ for all $K$, then $(\phi^{-1})^*F\in\o_I$. Changing the 
variables $z=\phi(w)$, one gets 
$I_K=\int_{\phi^{-1}(K+i\R)}|F(w)|\,|\mbox{Jac}_\phi(w)|\,|dw\wedge\overline{dw}|$ 
where $\mbox{Jac}_\phi$ is the Jacobian determinant of $\phi$. A computation 
yields $\mbox{Jac}_\phi(w)=\cos^2(x)+\sinh^2(y)\quad(w=x+iy)$ and, since 
$\phi^{-1}(K+i\R)\subset]-\frac{\pi}{2},
\frac{\pi}{2}[+i\R$, one gets $I_K\leq\int^{\frac{\pi}{2}}_{-\frac{\pi}{2}}
\int^{\infty}_{-\infty}\exp\left(-a|y|^{\frac{1}{\beta}}
+b|x|^{\frac{1}{1-\alpha}}\right)\,\left(\cos^2(x)+\cosh^2(y)\right)\,dy\,dx$ 
which is lower than 
$2b\exp\left[\left(\frac{\pi}{2}\right)^{\frac{1}{1-\alpha}}\right]\int^{\infty}_{-\infty}
\exp\left(-a|y|^{\frac{1}{\beta}}\right)\,(1+\cosh^2(y))\,dy$, thus finite 
as soon as $\beta<1$. The exact same argument works for 
$I\subset]-\frac{\pi}{2},0[$, hence $(\phi^{-1})^*F\in\o_I\cap\o_{-I}$. A 
similar argument yields item (ii).
\EPf
\pagebreak
\section{Twisted Weyl product}
\subsection{The product formula}
Let us consider the fundamental space 
$\s^{\sigma(\alpha_1,\alpha_2)}_{(\alpha_1,\alpha_2)}(2)\qquad(\alpha_1,\alpha_2)\in\R^{2}$ 
(cf. Section~\ref{TYPE-S}). 
Let $\varphi\in\sA$ and consider the partial function 
$\varphi_{x_1}:x\mapsto\varphi(x_1,x)$. For all $x_1\in\R$, the function 
$\varphi_{x_1}$ belongs to ${\cal S}^{\alpha_{1}}_{\alpha_{2}}(1)=:{\cal S}^{\alpha_{1}}_{\alpha_{2}}$. Therefore provided some 
restrictions on $(\alpha_1,\alpha_2)$, the function 
$\LT^{-1}(\phi^{-1}_{q,k\frac{\pi}{2}})^*\LT(\varphi_{x_1})$ ($k=0,1$) is 
well defined as an element of $\s'_I$ (cf. Proposition~\ref{ISOTWIST}).
\begin{dfn}
We define the linear map 
$$
\sA\stackrel{\tau^{(k)}_q}{\longrightarrow}C^\infty(\R^2)\quad(k=0,1)
$$
by $$\tau^{(k)}_q:=\mbox{\rm id}_{x_1}\otimes\left(\LT^{-1}\circ
(\phi^{-1}_{q,k\frac{\pi}{2}})^*\circ\LT\right)_{x_2}\qquad(x_1,x_2)\in\R^2.$$
We denote by $E^{(k)}_{(\alpha_1,\alpha_2)}$ its range in $C^\infty(\R^2)$. The 
inverse map $\mbox{\rm id}_{x_1}\otimes\left(\LT^{-1}\circ
\phi_{q,k\frac{\pi}{2}}^*\circ\LT\right)_{x_2}|_{\Ekab}$ will be denoted 
by $\Tkq$. It yields a linear isomorphism $\Tkq:\Ekab\to\sA$.
\end{dfn}
\begin{prop}
The formula 
$$
u\skq v:=\taukq\left(\Tkq u\sw\Tkq v\right)
$$
defines an associative $\R$-algebra structure $\skq$ on $\Ekab$.
\end{prop}
\Pf
One has $\Tkq u,\Tkq v\in\sA$ provided $u,v\in\Ekab$. We know that $\sA$ is 
stable under Weyl's product $\sw$ (cf. Proposition~\ref{WEYL}). Hence 
$u\skq v$ is well defined as an element of $\Ekab$. The associativity is 
obvious since $\skq$ is nothing else than the transportation of Weyl's
product via the isomorphism $\Tkq:\Ekab\to\sA$.
\EPf
\begin{dfn}
The product $\skq$ on $\Ekab$ will be referred as the twisted Weyl product.
\end{dfn}
\subsection{Observables of exponential type}
Given a convergent star product, an important question is the one of 
existence of non-tempered observables in the domain of the star product. 
This question has been studied for the case of the Moyal-Weyl product in 
\cite{OMMY}.

Let $A=(\alpha_1,\alpha_2)\in]0,1[^2$ with $\alpha_{1}+\alpha_{2}\geq1$. Set $\s_A:=\sA$. 
Viewing $\s_A=\left(\s_A\right)_x\quad(x\in\R^2)$ as a subspace of 
$C^\infty(\R^2)$, we consider the following alternate presentation of 
$\s_A$. Consider the sequence of maps 
$$
\begin{array}{ccccc}
\left(\s_A\right)_x& \to&\o(\C^2)&\to&C^\infty(\R^2)\\
f(x_1,x_2)&\mapsto&f(z_1=x_1+iy_1, z_2=x_2+iy_2)&\mapsto&f(iy_1,x_2).
\end{array}
$$
The function $\hat{f}(y_1,x_2):=f(iy_1,x_2)$ determines completely $f$. So 
that we have an injection $\left(\s_A\right)_x\to 
C^\infty(\R^2):f\mapsto\hat{f}$. Remark that the space $\hat{\s_A}$ 
contains elements of exponential growth. For example, one has
$f=e^{-(z_1^2+z_2^2)}\in\s_{(\frac{1}{2},\frac{1}{2})}$, which yields 
$|\hat{f}(y_1,x_2)|=e^{y_1^2-x_2^2}$.

We denote by $\widehat{\sw}$ the product on $\hat{\s_A}$ obtained by 
transporting Weyl's product on $\left(\s_A\right)_x$ via $f\mapsto\hat{f}$.
Observe that $\hat{\s_A}$ is still stable under the pointwise 
multiplication whose $\widehat{\sw}$ is a non-commutative deformation of. 
Observe also that for every $\psi\in\hat{\s_A}$ and $y_1\in\R$ the 
partial function $x_2\mapsto\psi(y_1,x_2)$ is in ${\cal S}^{\alpha_{1}}_{\alpha_{2}}$. Therefore the 
transformations $\taukq$ and $\Tkq$ are well defined on $\hat{\s_A}$.
\begin{prop}
Set $\widehat{\Ekab}:=\taukq(\hat{\s_A})$. Then for all $a,b\in\widehat{\Ekab}$, 
the formula 
$$
a\;\widehat{\skq}\;b:=\taukq\left(\Tkq a\;\widehat{\sw}\;\Tkq b\right)
$$
defines an associative $\R$-algebra structure on $\widehat{\Ekab}$. The space 
$\widehat{\Ekab}$ contains elements of exponential growth.
\end{prop}

\end{document}